\documentclass[12pt]{article}
\usepackage{amsmath,amsfonts,amssymb,amsthm,amscd}
\title{Stable embeddedness and $NIP$}
\date{January 4, 2010}
\author{Anand Pillay\thanks{Supported
by EPSRC grant EP/F009712/1}\\University of Leeds}
\newtheorem{Theorem}{Theorem}[section]
\newtheorem{Proposition}[Theorem]{Proposition}
\newtheorem{Definition}[Theorem]{Definition} 
\newtheorem{Remark}[Theorem]{Remark}
\newtheorem{Lemma}[Theorem]{Lemma}
\newtheorem{Corollary}[Theorem]{Corollary}
\newtheorem{Fact}[Theorem]{Fact}

\begin{document}
\maketitle

\begin{abstract} 
We give some sufficient conditions for a predicate $P$ in a complete theory $T$ to be ``stably embedded". Let $\cal P$ be $P$ with its ``induced $\emptyset$-definable structure". The conditions are that ${\cal P}$ (or rather its theory) has finite thorn rank, $P$ has $NIP$ in $T$ and that $P$ is stably $1$-embedded in $T$.  This generalizes a recent result of Hasson and Onshuus \cite{Hasson-Onshuus} which deals with the case where $P$ is $o$-minimal in $T$. Our proofs make use of the theory of strict nonforking and weight in $NIP$ theories  (\cite{Chernikov-Kaplan}, \cite{Shelah783}).
\end{abstract}

\section{Introduction and preliminaries}
The notion of ``stable embeddedness" of a predicate (or definable set) in a theory (or structure) is rather important in model theory, and says roughly that no new structure is added by external parameters.  The current paper is somewhat technical, and heavily influenced by \cite{Hasson-Onshuus} on the stable embeddability of $o$-minimal structures, which we wanted to cast in a more general context.

Stable embeddedness usually refers to a structure $M$ which is interpretable (without parameters) in another structure $N$, and says that any subset of $M^{n}$ which is definable, with parameters, in $N$, is definable, with parameters in $M$. I prefer to think of the universe of $M$ as simply an $\emptyset$-definable set $P$ in $N$, and to say that $P$ is stably embedded in $N$ if every subset of $P^{n}$ which is definable (in $N$) with parameters, is definable (in $N$) with parameters from $P$. See below.

So our general framework is that of a complete theory $T$ in language $L$ and a distinguished predicate or formula or sort $P(x)$.  We work in a saturated model $N$ of $T$, and
let $P$ also denote the interpretation of $P$ in $N$. Unless we say otherwise, ``definability" refers to definability with parameters in the ambient structure $N$. 
We will also discuss the notion ``$P$ is stable in $N$" just to clarify current notation and relationships.

\begin{Definition} (i) $P$ is stable in $T$ (or in $N$) if there do NOT exist a formula $\phi({\bar x},y)$  (where ${\bar x}$ is a tuple of variables each of which is of sort $P$, and $y$ is an arbitrary tuple of variables) and ${\bar a_{i}} \subset P$ and $b_{i} \in N^{eq}$ for $i<\omega$ such that $N\models \phi({\bar a_{i}},b_{j})$ iff $i\leq j$. 
\newline 
(ii)  $P$ is $NIP$ in $T$ (or $N$) if there do NOT exist $\phi({\bar x},y)$ (with same proviso as before)  and ${\bar a_{i}}\subset P$ for $i<\omega$ and $b_{s}\in N^{eq}$ for  $s\subseteq \omega$  such that $N\models \phi({\bar a_{i}},b_{s})$ iff $i\in s$. 
\newline
(iii) $P$ is stably embedded in $T$ (or $N$) if for all $n$ every subset of $P^{n}$ which is definable in $N$ with parameters, is definable in $N$ with parameters from $P$.
\newline
(iv) $P$ is $1$-stably embedded, if  (iii) holds for $n=1$.

\end{Definition}

We let ${\cal P}$ denote the structure whose universe is $P$ and whose basic relations are those which are $\emptyset$-definable in $N$.
Note that if $P$ is $NIP$ in $N$ and ${\cal P}$ is the structure with universe $P$ and relations all subsets of various $P^{n}$ which are $\emptyset$-definable in $N$, then $Th({\cal P})$ has $NIP$ too. Of course if $T$ (= $Th(N)$) has $NIP$ then $P$ has $NIP$ in $N$, and even in this case there are situations where $P$ is not stably embedded in $N$. For example when $T$ is the theory of dense pairs of real closed fields and where $P$ is the bottom model.

\begin{Remark}
(i) For $P$ to be stable in $N$ it is enough that Definition 1.1(i) holds in the case where ${\bar x}$ is a single variable $x$
ranging over $P$. Likewise for $P$ being $NIP$ in $N$ and Definition 1.1(ii). Also if $P$ is stable in $N$ then it is $NIP$ in $N$.
\newline
(ii) Suppose that $<$ is a distinguished $\emptyset$-definable total ordering on $P$. Define $P$ to be $o$-minimal in $N$ if every definable (in $N$, with parameters) subset of $P$ is a finite union of intervals (with endpoints in $P$ together with plus or minus $\infty$) and points. Then IF $P$ is $o$-minimal in $N$ then it is $1$-stably embedded in $N$ AND has $NIP$ in $N$. 
\newline
(iii) If $P$ is stable in $N$ then $P$ is stably embedded in $N$.
\end{Remark}
\noindent
{\em Comments.} (i) (for $P$ being $NIP$ in $T$) and (ii) were already mentioned in an earlier draft of \cite{Hasson-Onshuus}. For (i) the point is that well-known results 
reducing stability and $NIP$ to the case of formulas $\phi(x,y)$ where $x$ ranges over elements rather than tuples, adapt to the current ``relative" context. (iii) is also well-known.

\vspace{2mm}
\noindent
We now discuss the finite thorn rank condition. There is an extensive theory in place around rosy theories and $U$-thorn-rank (e.g. \cite{Onshuus}). 
From this point of view, a theory $T'$ has ``finite rank" if $T'$ is rosy and there is a finite bound on the $U$-thorn-ranks of types 
in any given imaginary sort. If $T'$ is $1$-sorted it suffices to have a bound on the $U$-thorn-ranks of types of elements in the home sort. See 
\cite{EKP} for more details on this point of view.
We will 
now give an equivalent abstract definition which is the one we will work with and which does not need any acquaintance with rosy theories and/or thorn forking.
In the following $T'$ is a complete theory,  $a,b,c,..$ range over possibly imaginary elements of a saturated model of $T'$ and $A,B,C,..$ range over small sets of imaginaries from such a model.
\begin{Definition} We say that $T'$ has {\em finite rank} if there is an assignement to every $a,B$ of a nonnegative integer $rk(a/B)$, depending only on
$tp(a,B)$ with the following properties:
\newline
(i) $rk(a,b/C) = rk(b,a/C) = rk(a/b,C) + rk(b/C)$,
\newline
(ii) $rk(a/B) = 0$ iff $a\in acl(B)$,
\newline
(iii) for any $a$ and $B\subseteq C$, there is $a'$ such that $tp(a'/B) = tp(a/B)$ and $rk(a'/C) = rk(a/B)$,
\newline
(iv) for any (imaginary) sort $Z$ there is a finite bound on $rk(a)$ for $a\in Z$. 
\newline

\end{Definition}

If $T'$ is a finite rank theory then we obtain a good notion of independence by defining $a$ to be independent from $B$ over $C$ if 
$rk(a/BC) = rk(a/C)$. Of course an $o$-minimal theory is an example.

\vspace{2mm}
\noindent
We can now state our main result, reverting to our earlier notation ($T$,$P$,${\cal P}$, etc.). 
\begin{Theorem} Suppose that $Th({\cal P})$ is finite rank, $P$ has $NIP$ in $T$ and $P$ is $1$-stably embedded in $N$. Then $P$ is stably embedded in $N$.
\end{Theorem}



\vspace{2mm}
\noindent
Our previous discussion shows that if $P$ is $o$-minimal in $T$ then all the hypotheses of Theorem 1.4 are satisfied.
Let us note immediately that Theorem 1.4 needs the ``$NIP$ in $N$" hypothesis on $P$. For example let $N$ be the structure with two disjoint unary predicates $P, Q$, and a random bi-partite graph relation $R$ between $P^{(2)}$ (unordered pairs from $P$) and $Q$. Then one checks that $P$ is $1$-stably embedded in $N$ but not stably embedded. 

\vspace{2mm}
\noindent
Hrushovski suggested the following example showing that the finite rank hypothesis is needed. Let $T$ be some completion of the theory of proper dense elementary pairs of algebraically closed valued fields  $F_{1} < F_{2}$, which have the same residue field and value group. Then $F_{1}$ is $1$-stably embedded in $F_{2}$. But
if $a \in F_{2}\setminus F_{1}$, then the function taking $x\in F_{1}$ to $v(x-a)$ cannot be definable with parameters from $F_{1}$.

\vspace{2mm}
\noindent
Thanks to Ehud Hrushovski and Alf Onshuus for helpful comments, suggestions, and discussions.

\section{Forking} 
In this section we fix a complete theory $T$ and work in a saturated model ${\bar M}$. As usual $A,B,..$ denote small subsets of ${\bar M}$. Likewise for small elementary substructures $M, M_{0}$ etc. $a,b$,.. usually denote
elements of ${\bar M}^{eq}$  unless we say otherwise. Likewise for variables. Dividing and forking are meant in the sense of Shelah. Namely a formula $\phi(x,b)$ divides over $A$ if some $A$-indiscernible sequence $(\phi(x,b_{i}):i<\omega)$ (with $b_{0} = b$) is inconsistent. And a 
partial type forks over $A$ if it implies a finite disjunction of formulas each of which divides over $A$. By a global type we mean a complete type over ${\bar M}$ (or over a sufficiently saturated $M$). Note that for a global type $p(x)$, $p$ does not fork over $A$ iff $p$ does not divide over $A$ (i.e. every formula in $p$ does not divide over $A$). Also any partial type does not fork over $A$ if and only if it has an extension to a global type which does 
not divide (fork) over $A$. Let $M_{0}$ denote a small model (elementary substructure of ${\bar M}$). 
\begin{Fact} (i) Suppose that $T$ has $NIP$. Let $p(x)$ be a global type. Then $p$ does not fork over $M_{0}$ if and only if $p$ is $Aut({\bar M}/M_{0})$-invariant.
\newline
(ii) If the global type $p(x)$ is $Aut({\bar M}/M_{0})$-invariant, and $(a_{i}:i<\omega)$ are such that $a_{n+1}$ realizes
$p|(M_{0}a_{0}..a_{n})$ for all $n$, then $(a_{i}:i<\omega)$ is indiscerinible over $M_{0}$ and its type over $M_{0}$ depends only on $p$.
$(a_{i}:i<\omega)$ is called a Morley sequence in $p$ over $M_{0}$.
\end{Fact}
\noindent
{\em Comment.}  (i) appears explicitly in \cite{Adler}, and also in \cite{NIPII}  but is implicit in 
\cite{Shelah783}.
\newline
(ii) is well-known, but see Chapter 12 of \cite{Poizat} for a nice account.

\begin{Definition} (i) Let $p(x)$ be a global type (or complete type over a saturated model), realized by $c$. We say that $p$ strictly does not fork over $A$ if $p$ does not fork over $A$ and $tp({\bar M}/Ac)$ does not fork over $A$  (namely for each small $B$ containing $A$ and realization $c'$ of $p|B$, $tp(B/Ac)$ does not fork over $A$).
\newline
(ii) Let $A\subseteq B$, and $p(x)\in S(B)$. Then $p$ strictly does not fork over $A$ (or $p$ is a strict nonforking extension of $p|A$) if $p$ has an extension to a global complete type $q(x)$ which strictly does not fork over $A$.
\end{Definition}

\begin{Fact} Assume that $T$ has $NIP$ and let $M_{0}$ be a model.
\newline
(i) For any formula $\phi(x,b)$, $\phi(x,b)$ divides over $M_{0}$ iff $\phi(x,b)$ forks over $M_{0}$.
\newline
(ii) Any $p(x)\in S(M_{0})$ has a global strict nonforking extension.
\newline
(iii) Suppose $p(x)$ is a global complete type which strictly does not fork over $M_{0}$. Let $(c_{i}:i<\omega)$ be a Morley sequence in $p$ 
over $M_{0}$. Suppose $\phi(x,c_{0})$ divides over $M_{0}$ (where $\phi(x,y)$ is over $M_{0})$. Then $\{\phi(x,c_{i}):i<\omega\}$ is inconsistent.
\end{Fact}
\noindent
{\em Proof.} This is all contained in \cite{Chernikov-Kaplan}, where essentially only the property $NTP_{2}$ (implied by $NIP$) is used. (i) is Theorem 1.2 there (as any model is an ``extension base" for nonforking). (ii) is Theorem 3.29 there. And (iii) is Claim 3.14 there.

\begin{Corollary} Suppose that $T$ has $NIP$ and $M_{0}$ is a model of cardinality $\kappa \geq |T|$. Let $q(y)$ be a global complete type which strictly does
not fork over $M_{0}$ and let $(c_{\alpha}:\alpha < \kappa^{+})$ be a Morley sequence in $q$ over $M_{0}$. Then for any (finite) tuple $a$, there is $\alpha<\kappa^{+}$ such that $tp(a/M_{0}c_{\alpha})$ does not fork over $M_{0}$.
\end{Corollary}
\noindent
{\em Proof.} Suppose not. So, using Fact 2.3(i), for each $\alpha$ there is some formula $\phi_{\alpha}(x,y)$ with parameters from $M_{0}$ such that 
$\models \phi_{\alpha}(a,c_{\alpha})$ and  $\phi_{\alpha}(x,c_{\alpha})$ divides over $M_{0}$. As $M_{0}$ has cardinality $\kappa\geq |T|$ there is $\phi(x,y)$ over $M_{0}$ and $\alpha_{0} < \alpha_{1} < \alpha_{2} < ...  < \kappa^{+}$ such that $\phi(x,y) = \phi_{\alpha_{i}}(x,y)$ for all 
$i< \omega$. By Fact 2.3 (iii) $\{\phi(x,c_{\alpha_{i}}):i<\omega\}$ is inconsistent, which is a contradiction as this set of formulas is supposed to be realized by $a$.

\vspace{5mm}
\noindent
Some remarks are in order concerning the notions introduced above and the last corollary. Strict nonforking in the $NIP$ context was introduced by Shelah \cite{Shelah783} (where the study of forking in $NIP$ theories was also initiated) and all the results above are closely connected in one way or another with Shelah's work. A version of Corollary 2.4, which on the face of it is incorrect due possibly to typographical errors, appears in \cite{Shelah783} as Claim 5.19. A better version of Corollary 2.4 appears in \cite{Usvyatsov}.

\section{Proof of Theorem 1.4}
We revert to the context of section 1, and Theorem 1.4.  Namely $T$ is an arbitrary theory, $N$ a saturated model, and $P$ a $\emptyset$-definable set in 
$N$.  As there ${\cal P}$ is the structure with universe $P$ and relations subsets of various $P^{n}$ which are $\emptyset$-definable in $N$.
So if for example  $a$ and $b$ are tuples from $P$ (or even from ${\cal P}^{eq}$) which have the same type in ${\cal P}$ then they have the same type in $N$.  We will assume that $Th({\cal P})$ is finite rank, as in Definition 1.3. We assume that $P$ has $NIP$ in $N$ (from which it follows that any sort in ${\cal P}^{eq}$ has $NIP$ in $N$ and also the structure ${\cal P}$ has $NIP$ in its own right).  We also assume that $P$ is $1$-stably embedded in $N$. Note that any sort (or definable set) in ${\cal P}^{eq}$ is also a sort (or definable set) in $N^{eq}$. If $S$ is a sort of ${\cal P}^{eq}$, and $X$ is a subset of $S$ which is definable (with parameters) in $N$, we will say, hopefully without ambiguity, that $X$ is {\em coded in $P$} if $X$ can be defined in $N$ with parameters from $P$, which is equivalent to saying that a canonical parameter for $X$ can be chosen in ${\cal P}^{eq}$, and is also of course equivalent to saying that $X$ is definable with parameters in the structure ${\cal P}$. Our aim is to prove that for any $n$, any subset of $P^{n}$ definable in $N$ is coded in $P$.
The $1$-stable embeddedness of $P$ states that any subset of $P$ definable with parameters in $N$ is coded in $P$. We first aim towards the following key proposition, which extends this to definable functions on $P$.

\begin{Proposition} Let $Z$ be a sort in ${\cal P}^{eq}$. Let $f:P\to Z$ be a function definable in $N$. Then $f$ is coded in $P$.
\end{Proposition}

\vspace{2mm}
\noindent
{\em Proof of Proposition 3.1.}  
\newline
This will go through a couple of steps. The first is an adaptation of ideas from \cite{Hasson-Onshuus}, using just the finite rank hypothesis on ${\cal P}$:
\begin{Lemma} (In the situation of Proposition 3.1) There is a relation $R(x,z)$ which is definable in $N$ with parameters from $P$, and some $k<\omega$, such that
\newline
$N\models  (\forall x\in P)(R(x,f(x)) \wedge (\exists ^{\leq k}z\in Z)R(x,z))$
\end{Lemma} 
\noindent
{\em Proof.}  Note that we are free to add parameters from $P$ to the language (but not, of course, from outside $P$). We try to construct  $a_{i}\in P$, $b_{i}\in Z$ for $i< \omega$ such that
\newline
(i) $a_{n} \notin acl((a_{i},b_{i})_{i<n})$,
\newline
(ii) $b_{n} = f(a_{n})$, and
\newline
(iii) $b_{n}\notin acl((a_{i}b_{i})_{i<n},a_{n})$.

\vspace{2mm}
\noindent
If at stage $n$ we cannot continue the construction it means that for all $a\in P$ such that $a\notin acl((a_{i}b_{i})_{i<n}$, 
$f(a)\in acl((a_{i}b_{i})_{i<n},a)$. Compactness will yield the required relation $R(x,z)$  (defined with parameters from $P$). 
\newline
So we assume that the construction can be carried out and  aim for a contradiction. We now work in the finite rank structure
${\cal P}$. Note that $rk(a_{n}b_{n}/(a_{i}b_{i})_{i<n}) > 1$ for all $n$ and moreover there is a finite bound to $rk(a_{n}b_{n})$ as $n$ varies.
Choose minimal $s> 1$ such that for some $m$, for infinitely many $n > m$, $rk(a_{n}b_{n}/(a_{i}b_{i})_{i\leq m}) = s$. It is then easy to find
$m < i_{0} < i_{i} < .. $ such that for all $j$, 
\newline
$rk(a_{i_{j}},b_{i_{j}}/a_{0},b_{0},..,a_{m},b_{m},a_{i_{0}},b_{i_{0}},..,a_{i_{j-1}},b_{i_{j-1}}) = 
rk(a_{i_{j}},b_{i_{j}}/a_{0},b_{0},..,a_{m},b_{m}) = s$.
\newline
So after adding constants (in $P$) for $a_{0},b_{0},..,a_{m},b_{m}$, thinning the sequence, and relabelling, we have 
\newline
(iv) $rk(a_{n},b_{n}) = rk(a_{n}b_{n}/(a_{i}b_{i})_{i<n}) = s > 1$ for all $n$, which means that
\newline
(v) $\{(a_{0}b_{0}), (a_{1}b_{1}), (a_{2}b_{2}),....\}$ is an ``independent" set of tuples.

\vspace{2mm}
\noindent
Of course we still have that $a_{i}\notin acl(\emptyset)$, $b_{i}\notin acl(a_{i})$ for all $i$, as well as $b_{i} = f(a_{i})$
for all $i$.  Now $f$ is definable in $N$ with some parameter $e$ so let us write $f$ as $f_{e}$ to exhibit the dependence on $e$.
\newline
{\em Claim.}  For each $S\subseteq \omega$  there are $b_{i}^{S}\in Z$ for $i<\omega$ such that
\newline
(a) For $i\in \omega$, $i\in S$ iff $b_{i}^{S} = b_{i}$.
\newline
(b)  $tp((a_{i}b_{i})_{i<\omega}) = tp((a_{i}b_{i}^{S})_{i<\omega})$.
\newline
{\em Proof of claim.}  Fix $n\notin S$. By (v) above, $rk(b_{n}/a_{n},(a_{i}b_{i})_{i\neq n}) > 0$, so let $b_{n}^{S}$ be such that 
$tp(b_{n}^{S}/a_{n},(a_{i}b_{i})_{i\neq n}) =  tp(b_{n}/a_{n},(a_{i}b_{i})_{i\neq n})$ and $rk(b_{n}^{S}/(a_{i}b_{i})_{i<\omega}) > 0$. So we see that $tp((a_{n}b_{n}),(a_{i}b_{i})_{i\neq n}) = tp(a_{n}b_{n}^{S}),(a_{i}b_{i})_{i\neq n})$ and $b_{n}^{S} \neq b_{n}$. Iterate this to obtain the claim.

\vspace{2mm}
\noindent
By (b) of the claim, and automorphism, for each $S\subseteq \omega$, there is $e_{S}$ such that  $f_{e_{S}}(a_{n}) = b_{n}^{S}$ for all $n\in \omega$. 
By (a) of the claim, we have that $f_{e_{S}}(a_{n}) = b_{n}$ iff $n\in S$. This contradicts $P$ having $NIP$ in $T$, and proves Lemma 3.2

\vspace{5mm}
\noindent
Note that if $Th({\cal P})$ had Skolem functions, or even Skolem functions for ``algebraic" formulas, we could quickly deduce  Proposition 3.1 from Lemma 3.2. (And this is how it works in the $o$-minimal case.)
Likewise if we can choose $R$ in Lemma 3.2 such that $k=1$, we would be finished. So let us choose $R$ in Lemma 3.2 such that $k$ is minimized, and work towards showing that $k = 1$. Here we use nonforking in the $NIP$ theory $Th({\cal P})$.  Let us fix a small elementary substructure $M_{0}$ of ${\cal P}$ which contains the parameters from $R$. So for any $a\in P$, $f(a)\in acl(M_{0},a)$, and this will be used all the time.

\begin{Lemma} There is a finite tuple $c$ from ${\cal P}^{eq}$ such that whenever $a\in P$ and $b = f(a)$ and $tp(a/M_{0}c)$ does not fork over $M_{0}$ (in the structure ${\cal P}$) then $b\in dcl(M_{0},c,a)$  (in ${\cal P}$ or equivalently in $N$).
\end{Lemma}
\noindent
{\em Proof.}  Let us suppose the lemma fails (and aim for a contradiction). We will first find inductively $a_{n},b_{n},d_{n}$ for $n< \omega$ such that writing $c_{n} = (a_{n},b_{n},d_{n})$ we have
\newline
(i) $a_{n}\in P$, $b_{n}\neq d_{n}$ are in $Z$ and $b_{n} = f(a_{n})$.
\newline
(ii) $tp(a_{n}b_{n}/M_{0}c_{0}..c_{n-1})$ does not fork over $M_{0}$, and
\newline
(iii) $tp(a_{n}b_{n}/M_{0}c_{0}...c_{n-1}) = tp(a_{n}d_{n}/M_{0}c_{0}...c_{n-1})$.

\vspace{2mm}
\noindent
Suppose have found $a_{i},b_{i},d_{i}$ for $i<n$. Put $c = c_{0}...c_{n-1}$  (so if $n=0$, $c$ is the empty tuple.) As the claim fails for $c$, we find $a\in P$, such that $f(a) = b$, $tp(a/M_{0}c)$ does not fork over $M_{0}$, and $b\notin dcl(M_{0}ca)$. Noting that $b\in acl(M_{0},a)$ we also have that $tp(ab/M_{0}c)$ does not fork over $M_{0}$. Let $d$ realize the same type as $b$ over $M_{0}ca$ with $d\neq b$. Put $a_{n} = a$, $b_{n} = b$, $d_{n} = d$. So the construction can be accomplished.

 \vspace{2mm}
 \noindent
 {\em Claim.} Let $S \subseteq \omega$. Define $g_{n} = b_{n}$ if $n\in S$ and $g_{n} = d_{n}$ if $n\notin S$. Then $tp((a_{i}g_{i})_{i<\omega}/M_{0}) = tp((a_{i}b_{i})_{i<\omega}/M_{0}$ (in ${\cal P}$ or equivalently in $N$).
 \newline
 {\em Proof of claim.} Assume, by induction, that 
\newline
(*) $tp((a_{i}g_{i})_{i<n}/M_{0}) = tp((a_{i}b_{i})_{i<n}/M_{0})$, 
\newline
and we want the same thing with $n+1$ in place of $n$. By (iii) above it suffices to prove that $tp((a_{i}g_{i})_{i<n}(a_{n}b_{n})/M_{0}) = tp((a_{i}b_{i})_{i<n}(a_{n}b_{n})/M_{0})$. This is true by (*) and Fact 2.1(i). (Namely as $tp(a_{n}b_{n}/M_{0}c)$ does not fork over $M_{0}$ it has a global complete extension $q$ say which does not fork over $M_{0}$ so by 2.1(i) $q$ is fixed by automorphisms which fix $M_{0}$ pointwise. So whether or not a formula $\phi(x,z,d)$ say is in $q$ or not depends just on $tp(d/M_{0})$, so the same is true for 
$tp(a_{n}b_{n}/M_{0}c)$.) The claim is proved.

\vspace{2mm}
\noindent
As earlier we write $f$ as $f_{e}$ where $e$ is a parameter in $N$ over which $f$ is defined.
By the claim, and automorphism, for each $S\subseteq \omega$ we can find $e_{S}$ in $N$ such that $f_{e_{S}}(a_{i}) = g_{i}$ for all $i< \omega$. In particular, as $d_{i}\neq b_{i}$ for all $i$,  $f_{e_{S}}(a_{i}) = b_{i}$ if and only if $i\in S$, showing that $P$ has the independence property in $N$, a contradiction.
Lemma 3.3 is proved.

\vspace{5mm}
\noindent
We now complete the proof of Proposition 3.1 by showing that $k=1$. Let us assume $k > 1$ and get a contradiction. 
Let $c\in {\cal P}^{eq}$ be as given by the Lemma 3.3. Then, (by the previous two lemmas), whenever $a\in P$ and $tp(a/M_{0}c)$ does not fork over $M_{0}$, THEN there is 
$d\in Z$ such that $\models R(a,d)$ and $d\in dcl(M_{0},a,c)$.  Note that $f$ is not mentioned here, so this statement is purely about ${\cal P}$.
Hence:
\newline
(*) whenever $c'\in {\cal P}$ has the same type over $M_{0}$ as $c$, then whenever $a\in P$ and $tp(a/M_{0}c')$ does not fork over $M_{0}$ then 
$d\in dcl(M_{0},a,c)$ for some $d\in Z$ such that $\models R(a,d)$.

\vspace{2mm}
\noindent
We work in ${\cal P}$. By Fact 2.3 (ii) let  $q(y)$ be a global strict nonforking extension of $tp(c/M_{0})$, and let 
$(c_{\alpha}:i< |T|^{+})$ be a Morley sequence in $q$ over $M_{0}$. 
By Corollary 2.4, for any $a\in P$, there is $\alpha$ such that $tp(a/M_{0}c_{\alpha})$ does not fork over $M_{0}$, and hence by (*), there is $d$ such that $\models R(a,d)$ and $d\in dcl(M_{0},a,c_{\alpha})$. By compactness there is a partial function $g(-)$ defined over 
$M_{0}\cup \{c_{\alpha}:\alpha < |T|^{+}\}$ such that
\newline
(**) for all $a\in P$, there is $d\in Z$ such that $\models R(a,d)$ and $g(a) = d$.

\vspace{2mm}
\noindent
Now by our assumption that $P$ is $1$-stable embedded, $\{a\in P: g(a) = f(a)\}$ is definable over parameters from $P$ by some formula $\psi(x)$ say.
Let $R'(x,z)$ be   $(\psi(x)\wedge z= g(x)) \vee (\neg\psi(x)\wedge R(x,z)\wedge z \neq g(x))$. Then clearly $R'$ satisfies Lemma 3.2 with $k-1$. This contradiction proves Proposition 3.1.

\vspace{5mm}
\noindent
{\em Proof of Theorem 1.4.} 
\newline
There is no harm in adding a few constants for elements of $P$ (so that we can do definition by cases).
We will prove by induction on $n$ that any subset of $P^{n}$ which is definable in $N$ with parameters, is coded in ${\cal P}^{eq}$. For $n=1$ this is the $1$-stable embeddednes hypothesis. Assume true for $n$ and we'll prove for $n+1$.
Let $\phi(x_{1},..,x_{n},x_{n+1},e)$ be a formula, with parameter $e\in N$ defining the set $X\subseteq P^{n+1}$. For each $a\in P$, let $X_{a}$ be the subset of $P^{n}$ defined by
$\phi(x_{1},x_{2},..,x_{n},a,e)$. By induction hypothesis $X_{a}$ is coded in $P$, namely has a canonical parameter $c_{a}$ say in ${\cal P}^{eq}$. By compactness, we assume that there is a formula $\psi(x_{1},..,x_{n},z)$  (of $L$), where $z$ ranges over a sort $Z$ of ${\cal P}^{eq}$, such that for each $a\in P$, $X_{a}$ is defined by $\psi(x_{1},..,x_{n},c_{a})$  (and $c_{a}$ is still a canonical parameter for $X_{a}$). The map taking $a$ to $c_{a}$ is clearly an $e$-definable function $f$ say from $P$ to $Z$. By Proposition 3.1, $f$ is coded in ${\cal P}$  (i.e. can be defined with parameters from $P$). As the original
set $X\subseteq P^{n+1}$ is defined by the ``formula" 
$\psi(x_{1},..,x_{n},f(x_{n+1}))$, it follows that $X$ is also coded in $P$.  The proof of Theorem 1.4 is complete.



\begin{thebibliography}{99}
\bibitem{Adler} H. Adler, Introduction to theories without the independence property, to appear in Archive Math. Logic.


\bibitem{Chernikov-Kaplan} A. Chernikov and I. Kaplan, Forking and dividing in $NTP_{2}$ theories, to appear in Journal of Symbolic Logic.
\newline
(number 147 on the MODNET preprint server)

\bibitem{EKP}  C. Ealy, K. Krupinski, and A. Pillay, Superrosy dependent groups having finitely satisfiable generics, Annals of Pure and Applied Logic 151 (2008), 1 - 21.

\bibitem{Hasson-Onshuus}  A. Hasson and A. Onshuus, Embedded $o$-minimal structures, to appear in Journal of the  London Math. Soc.


\bibitem{NIPII} E. Hrushovski and A. Pillay, On $NIP$ and invariant measures, preprint 2009 (revised version). 
\newline
http://arxiv.org/abs/0710.2330

\bibitem{Onshuus}  A. Onshuus, Properties and consequences of thorn independence, Journal of Symbolic Logic, 71 (2006), 1-21.

\bibitem{Poizat} B. Poizat, A course in Model theory; an introduction to contemporary mathematical logic, Springer 2000.


\bibitem{Shelah783} S. Shelah, Dependent first order theories, continued, Israel J. Math. 173 (2009), 1-60.
\newline
http://arxiv.org/abs/0406440

\bibitem{Usvyatsov} A. Usvyatsov, Morley sequences in dependent theories, preprint.
\newline
http://arxiv.org/abs/0810.0733



\end{thebibliography}
\end{document}